\documentclass[submission]{FPSAC2025}

\newtheorem{thm}{Theorem}[section]

\newtheorem{defn}[thm]{Definition}

\usepackage{tikz}
\usepackage{adjustbox}

\newcommand*\linear[1]{\trimbox{-1 -.5 -1 -10}{\tikz[baseline=(char.base)]{

            \begin{scope}[ color                 = green]
            \draw [<->] (-2.2,0.0) -- (-1.2,0.0);
            \node[] at (-1.7,0.0) (char) {};
            \draw[red,fill=red] (-2.243301270189222,-0.024999999999999967) circle (0.05) (char) {};
            \draw[red,fill=red] (-1.156698729810778,-0.024999999999999967) circle (0.05);
            \draw [<->] (-3.3,0.0) -- (-2.3,0.0);
            \draw[red,fill=red] (-3.3433012701892215,-0.024999999999999967) circle (0.05);
            \draw[red,fill=red] (-2.256698729810778,-0.024999999999999967) circle (0.05);
            \draw [<->] (-4.4,0.0) -- (-3.4,0.0);
            \draw[red,fill=red] (-4.443301270189222,-0.024999999999999967) circle (0.05);
            \draw[red,fill=red] (-3.3566987298107778,-0.024999999999999967) circle (0.05);
            \draw [<->] (-5.5,0.0) -- (-4.5,0.0);
            \draw[red,fill=red] (-5.543301270189222,-0.024999999999999967) circle (0.05);
            \draw[red,fill=red] (-4.456698729810778,-0.024999999999999967) circle (0.05);
            \draw [<->] (-6.6,0.0) -- (-5.6,0.0);
            \draw[red,fill=red] (-6.643301270189221,-0.024999999999999967) circle (0.05);
            \draw[red,fill=red] (-5.556698729810778,-0.024999999999999967) circle (0.05);
            \draw [<->] (-7.699999999999999,0.0) -- (-6.699999999999999,0.0);
            \draw[red,fill=red] (-7.743301270189221,-0.024999999999999967) circle (0.05);
            \draw[red,fill=red] (-6.656698729810778,-0.024999999999999967) circle (0.05);

            \end{scope}
 
             \begin{scope}[ color                 = black
                 , execute at begin node = $\displaystyle
                 , execute at end node   = $]
                  \node at (-0.6,-0.4) {(0,6)} ;
                  \node at (-7.75,-0.4) {(6,0)} ;
                  \node at (-5.5,0.2526279441628825) {10020000} ;
                  \node at (-2.2,-.22) {01000002} ;
            \end{scope}

            }}}

\usepackage{lipsum}

\title[Power series for ASEP]{An explicit power series result for the two type ASEP}

\author[D. Ash]{David W. Ash\thanks{\href{mailto:dash@alumni.stanford.edu}{dash@alumni.stanford.edu}.
Thanks to Real Time Agents Inc of Pleasant Hill, CA, USA for partially funding this work.}\addressmark{1}}

\address{\addressmark{1}Real Time Agents, Inc. }

\received{\today}

\abstract{Research in combinatorics has often focused on the ASEP--asymmetric simple exclusion process. The
ASEP is inspired by processes in statistical mechanics, and involves particles of various species moving
around a lattice. The particles do not change species. In the present paper, based on earlier results of Mortimer and Prellberg
and others, we obtain
a new power series results for the two type ASEP.}


\usepackage[backend=bibtex]{biblatex}
\begin{document}

\maketitle

\section{Introduction}

In this paper, we provide an explicit power series formula for some specific cases of the two type periodic ASEP--the ASEP on a ring. We do so by building on earlier work
of Mortimer and Prellberg \cite{mortimer} and Elizalde \cite{elizalde}.
The ASEP (asymmetric simple exclusion process),
is a structure that has frequently been referred to in the combinatorics
literature. The ASEP consists of a one dimensional periodic lattice, with each point on the
lattice being populated with either a particle or a hole. At random intervals, each particle attempts to move either clockwise
or counterclockwise with different but fixed probabilites (hence the term `asymmetric'). Our focus will be on ASEP on a one
dimensional periodic lattice, which can also be thought of as the ASEP on a ring,
a form of Markov process as noted in \cite{multiline} by Corteel, Mandelshtam and Williams.

\section{Definitions}

Following \cite{multiline}, a partition $ \lambda $ and then the set $S(\lambda )$ all permutations of $\lambda$ 
can be defined. We can then define the multispecies asymmetric simple exclusion
process ASEP($\lambda$) as a Markov process on $S(\lambda )$ with certain specific transition probabilities:

\begin{defn}
For all partitions $\lambda$,
ASEP($\lambda$) is a Markov process on $S(\lambda )$. We let $t$ be a constant with $0\le t\le 1$. Suppose we have two 
permutations $\mu\in S(\lambda )$ and $\nu\in S(\lambda )$.
The transition probability, $P_{\mu,\nu}$, is then given by:

\begin{itemize}
  \item If $\mu=AijB$ and
$\nu=AjiB$ with $i\ne j$, then $P_{\mu,\nu}=\frac{t}{n}$ if
$i>j$ and $P_{\mu,\nu}=\frac{1}{n}$ if $j>i$.
  \item If $\mu=iAj$ and
$\nu=jAi$ with $i\ne j$, then $P_{\mu,\nu}=\frac{t}{n}$ if
$j>i$ and $P_{\mu,\nu}=\frac{1}{n}$ if $i>j$.
  \item If neither of the above conditions apply but $\nu\ne\mu$ then $P_{\mu,\nu}=0$. If $\nu=\mu$ then 
$P_{\mu,\mu}=1-\sum_{\nu\ne\mu}P_{\mu,\nu}$.
\end{itemize}
\end{defn}

Lattice paths within the ASEP have been previously explored by Elizalde in \cite{elizalde}. Elizalde develops a bijection between the
ASEP and lattice walks in simplicial regions. Elizalde focuses on bijections involving the one-type ASEP. We will extend
this to the two type of ASEP. Elizalde defines a simplicial region as an ordered tuple of nonnegative integers summing to a
constant. Lattice paths within such a region are defined, and these are shown to have a bijection with 
ASEP$(1,1,...,1,0,0,...,0)$ where there are
$d$ $1$'s and $L$ $0$'s. We would like to extend Elizalde's results to the multitype ASEP. Elizalde uses his bijection
to find the number
of paths in the ASEP, drawing upon earlier work by Mortimer and Prellberg \cite{mortimer}. Mortimer and Prellberg focus on a
triangular path
problem before switching to the simplified one-dimensional simplex, a line as shown in Figure 1.

\vspace{0in}
\begin{center}
\linear{0}
\linebreak
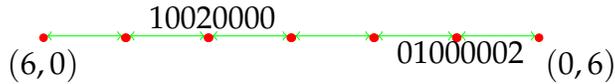
\textbf{Figure 1}: Simplified linear version of the simplex.
\end{center}

\section{Generating Functions}

To compute the generating function for certain cases of the one type ASEP, we can start with the following theorem of Mortimer and
Prellberg \cite{mortimer}, and later build up to the two type ASEP:

\begin{thm} \label{conj:name1}
The generating function $G_{u,v}(x,y;t)$ which counts $n$-step walks starting at $(u,v)$ in the simplified linear lattice
shown in Figure 1 is given by

\[ G_{u,v}(x,y;t) = \frac{1}{1-\frac{\frac{x}{y}+\frac{y}{x}}{p+\frac{1}{p}}} \left(x^uy^v 
                  - \frac{x^{u+v+1}p^{v+1}(1-p^{2u+2})}{y(1-p^{2u+2v+4})}
                  - \frac{y^{u+v+1}p^{u+1}(1-p^{2v+2})}{x(1-p^{2u+2v+4})}\right) \]

where $p$ is given by

\[ p = \frac{1-\sqrt{1-4t^2}}{2t} \]
\end{thm}

Here $x$ and $y$ count the position in the lattice
and $t$ counts the number of steps in the path. We can use this result to compute the generating function $G_{m,0,m,0}(1,1;t)$
which counts paths starting and ending at $(m,0)$. We observe first that $G_{m,0}(x,y;t)$ gives us the paths starting at
$(m,0)$ and ending anywhere. Such paths may be divided into two cases. The first case is that the path also ends at $(m,0)$.
Suppose the path returns to $(m,0)$ $k$ times, for $k\ge 0$, strictly between the start and the end. Consider the subpath between
two consecutive times when the path reaches $(m,0)$. The first step in such a path must be to $(m-1,1)$ and the last step must
be from $(m-1,1)$. Between the steps to and from $(m-1,1)$ the path must remain above $(m,0)$ implying that the generating
function for this part of the path is $G_{m-1,0,m-1,0}(1,1;t)$. As $k\ge 0$ this implies that the overall generating function
for this case is

\[ \sum_{k=0}^\infty G_{m-1,0,m-1,0}(1,1;t)^k \]

To label this with the fact that the path ends at $(m,0)$ this should be multiplied by $x^m$. The other subcase is that the path
ends somewhere other that $(m,0)$. In this subcase the above generating function covers up to the final appearance at $(m,0)$.
From there the path moves to $(m-1,1)$ and the remainder of the path is represented by the generating function
$G_{m-1,0}(x,y;t)$. To adjust for the fact that the original lattice has $u+v=m$ not $m-1$ we need to multiply by $y$. All this
gives us

\[ G_{m,0}(x,y;t) = (x^m+yG_{m-1,0}(x,y;t))\sum_{k=0}^\infty G_{m-1,0,m-1,0}(1,1;t)^k \]

or

\[ G_{m,0}(x,y;t) = \frac{x^m+yG_{m-1,0}(x,y;t)}{1-G_{m-1,0,m-1,0}(1,1;t)} \]
 
After replacing $m$ with $m+1$ and rearranging terms we have proven:

\begin{thm} \label{conj:name1}
The generating function $G_{m,0,m,0}(1,1;t)$ which counts $n$-step walks starting and ending at $(m,0)$ in the simplified linear lattice
shown in Figure 4 is given by

\[ G_{m,0,m,0}(1,1;t) = 1-\frac{x^{m+1}+yG_{m,0}(x,y;t)}{G_{m+1,0}(x,y;t)}  \]

\end{thm}

We next look at the generating function $G_{m,0,0,m}(1,1;t)$ counting paths starting at $(m,0)$ and ending at $(0,m)$. As before
we start by observing that $G_{m,0}(x,y;t)$ gives us the paths starting at $(m,0)$ and ending anywhere. We define a segment, of
which there will be zero or more in the path, as follows. The first segment is from the start of the path to the first occurrence
of $(0,n)$ if any. The second segment is from the end of the first segment to the next occurrence of $(n,0)$ if any. This
continues until the end of the path where there may be a section of the path not part of any segment. Suppose there are $k$
segments. Then the generating function for the part of the path in the segments, similar to above, will be given by
$G_{m-1,0,0,m-1}(1,1;t)^k$. There are then two cases to consider. If $k$ is even, then the last part of the path, after the final
segment, will have the generating function $G_{m-1,0}(x,y;t)$. To adjust for the fact that the original lattice has $u+v=m$ not $m-1$
we need to multiply by $x$. If $k$ is odd, then the last section of the path starts from the opposite end, giving the
generating function $G_{0,m-1}(x,y;t)$ which must be multiplied by $y$. All this gives us

\[ 
\begin{aligned}
G_{m,0}(x,y;t) = & (xG_{m-1,0}(x,y;t)+yG_{0,m-1}(x,y;t)G_{m-1,0,0,m-1}(1,1;t)) \\
                 & \cdot \sum_{k=0}^\infty G_{m-1,0,0,m-1}(1,1;t)^{2k} 
\end{aligned}
\]

or

\[ G_{m,0}(x,y;t) = \frac{xG_{m-1,0}(x,y;t)+yG_{0,m-1}(x,y;t)G_{m-1,0,0,m-1}(1,1;t)}{1-G_{m-1,0,0,m-1}(1,1;t)^{2}} \]

This expands to give us

\[ 
\begin{aligned}
G_{m,0}(x,y;t)G_{m-1,0,0,m-1}(1,1;t)^{2} & +yG_{0,m-1}(x,y;t)G_{m-1,0,0,m-1}(1,1;t) \\
                                         & +xG_{m-1,0}(x,y;t)-G_{m,0}(x,y;t)=0 
\end{aligned}
\]

After replacing $m$ with $m+1$ and using the quadratic formula we have proven:

\begin{thm} \label{conj:name1}
The generating function $G_{m,0,0,m}(1,1;t)$ which counts $n$-step walks starting at $(m,0)$ and ending at $(0,m)$ in the 
simplified linear lattice shown in Figure 4 is given by

\[ G_{m,0,0,m}(1,1;t) = \frac{
\sqrt{
\begin{aligned}
& y^2G_{0,m}(x,y;t)^2 \\
& -4G_{m+1,0}(x,y;t) \\
& \cdot (xG_{m,0}(x,y;t)-G_{m+1,0}(x,y;t))
\end{aligned}
}
-yG_{0,m}(x,y;t)}{2G_{m+1,0}(x,y;t)}  \]

\end{thm}

\section{Power Series for Two Type ASEP}

The results of the previous section, combined with the Elizalde bijection, gives us
a power series on the one type ASEP$(1,1,0,0,...,0)$ where there are $L$ $0$'s. We use those results to
establish a power series for the two type ASEP$(2,1,0,0,...,0)$ using a similar argument. Consider
paths starting at $(m,0)$. The behavior will be similar to that of the one type ASEP until the first time
when the $2$ and the $1$ switch places. As this can only happen when the path returns to $(m,0)$ this part of the
overall path will be governed by $G_{m,0,m,0}(1,1;t)$. After the switch in places, which we can mark with the
indeterminate $d$, the path needs to go to $(0,m)$ before an opportunity to switch places again occurs. This
assumes we only allow $1$ and $2$ to switch places in one direction. As the new switch in places reverses the
first one, this switch can be marked with $\frac{1}{d}$ and will be governed by the power series 
$G_{m,0,0,m}(1,1;t)$. The final segment, after the last switch in places, will be governed either by
$G_{m,0}(x,y;t)$ or $G_{0,m}(x,y;t)$.

Let $\lambda=(2,1,0,0,...,0)$ with $L=m$ zeroes. Denote the power series for ASEP$(\lambda)$ marking switches 
of the $1$ and $2$ with a $d$ by $A_\lambda (x,y;t;d)$. The above argument shows

\[ 
\begin{aligned}
A_\lambda (x,y;t;d) = G_{m,0}(x,y;t) + G_{m,0,m,0}(1,1;t) & (dG_{m,0}(x,y;t)+G_{m,0,0,m}(1,1;t)G_{0,m}(x,y;t)) \\
                                                          & \cdot \sum_{k=0}^\infty (G_{m,0,0,m}(1,1;t))^{2k} 
\end{aligned}
\]

Summing the infinite series gives us our main ASEP theorem:

\begin{thm} \label{conj:nameasep}
The power series $A_\lambda (x,y;t;d)$ counting $n$-step walks in ASEP$(\lambda)$ where
$\lambda=(2,1,0,0,...,0)$ with $m$ zeroes, starting with $1$, $2$
adjacent, is generated by

\[ A_\lambda (x,y;t;d) = G_{m,0}(x,y;t) + \frac{G_{m,0,m,0}(1,1;t) (dG_{m,0}(x,y;t)+G_{m,0,0,m}(1,1;t)G_{0,m}(x,y;t))}
                                                {1 - G_{m,0,0,m}(1,1;t)^2}   \]

\end{thm}

\section{Future work}

We would like to be able to generalize Theorem \ref{conj:nameasep} to cases with
more than two particles. This will requires a generalization of the results of \cite{mortimer} to give
a complete solution to higher order simplices beyond a straight line. As of now \cite{mortimer} offers
only partial results even for a triangle.

\printbibliography

\end{document}